\documentclass[12pt]{article}
\usepackage[dvips]{graphicx}
\usepackage{amsmath,amsthm,amsfonts,amssymb,amscd,epic,eepic}

\newtheorem{thm}{Theorem}[section]

\newtheorem{defn}[thm]{Definition}

\newtheorem{lem}[thm]{Lemma}
\newtheorem{prop}[thm]{Proposition}

\newtheorem{rem}[thm]{Remark}

\newtheorem{asmp}[thm]{Assumption}

\title{Orbifold aspects of the Longo-Rehren subfactors}

\author{
{\sc Nobuya Sato}\footnote{Supported by the
Grants-in-Aid for Scientific Research, JSPS.}\\
Department of Mathematics, Rikkyo University\\ 
Nishi-Ikebukuro, Tokyo, 171-8501, JAPAN\\
e-mail: {\tt nobuya@rkmath.rikkyo.ac.jp}}

\begin{document}

\maketitle

\begin{abstract}
In this article, we will prove that the subsectors of $\alpha$-induced 
sectors for $M \rtimes \hat{G} \supset M$ forms a modular category, where 
$M \rtimes \hat{G}$ is the crossed product of $M$ by the group dual 
$\hat{G}$ of a finite group $G$. In fact, we will prove that it is equivalent 
to M\"uger's crossed product. 
By using this identification, we will exhibit an orbifold aspect of the 
quantum double of $\Delta$(not necessarily non-degenerate) obtained from 
a Longo-Rehren inclusion $A \supset B_\Delta$ under certain assumptions. 

We will apply the above description of the quantum double of $\Delta$ to the 
Reshetikhin-Turaev topological invariant of closed 3-manifolds, and we obtain 
a simpler formula, which is a degenerate version of Turaev's theorem that 
the Reshetikhin-Turaev invariant for the quantum double of a modular 
category $\hat{\Delta}$ is the product of Reshetikhin-Turaev invariant of 
$\hat{\Delta}$ and its complex conjugate.
\end{abstract}

\section{Introduction}
Orbifold phenomena have sometimes appeared in subfactor theory. The 
first appearance was to construct subfactors with Dynkin diagrams of 
type $D_{2n}$ out of subfactors with Dynkin diagrams of type $A_{4n-3}$ 
through Dynkin diagram automorphisms \cite{Kawahigashi}. This method in 
consequence suggested the orbifold construction removes the degeneracy 
of braiding, and it is in fact proved in \cite{EK2}.  
Along the same line, in \cite{EK0} Evans and Kawahigashi extended the 
method of orbifold construction to the Hecke algebra subfactors of Wenzl 
\cite{Wenzl}. In a more sophisticated way, 
Goto defined an orbifold subfactor as a simultaneous crossed product by 
non-strongly outer automorphism with the trivial Loi invariant \cite{Goto}. 

There is another known way to remove degeneracy of braiding. That is the 
quantum double construction. Originally, it is the way to construct a 
higher symmetric Hopf algebra out of the initial Hopf algebra and its 
dual Hopf algebra. In subfactor context, we have Ocneanu's asymptotic 
inclusion $M_\infty \supset M \vee M^{op}$ constructed from a hyperfinite 
II$_1$ subfactor $N \subset M$. The finite system of $M_\infty$-$M_\infty$ 
bimodules obtained from an asymptotic inclusion is known to form a
modular category. In fact, it is the so called 
center construction in category theory \cite{Kassel}. However, this 
correspondence is not obvious because Ocneanu's construction of the finite system of 
$M_\infty$-$M_\infty$ bimodules makes an ingenious use of topological 
quantum field theory in three dimensions in the sense of Atiyah \cite{Atiyah}. 
For infinite factors, Longo and Rehren introduced an interesting 
subfactor nowadays called the Longo-Rehren subfactor. The Longo-Rehren 
subfactor produces the same tensor category as the one 
the asymptotic inclusion does \cite{Masuda1}. 

In his paper \cite{Izumi1}, Izumi examined and clarified the structure of 
a Longo-Rehren subfactor and its quantum double in a completely algebraic 
way, i.e., without using any help of TQFT. Moreover, he proved that 
the quantum double obtained from a Longo-Rehren inclusion is a modular 
category and further gave the description of modular $S$- and 
$T$-matrices in the language of sectors.
Thus, the quantum double in subfactors also provides a machinery 
to remove degeneracy of braiding. See \cite{EK2}, \cite{Masuda2} for the 
relationship between the orbifold construction and the asymptotic (or Longo-Rehren) 
inclusion.

In \cite{EK2}, Evans and Kawahigashi proved that the quantum double of 
a finite system of bimodules with non-degenerate braiding $\hat{\Delta}$ 
is equivalent to $\hat{\Delta} \otimes \hat{\Delta}^{op}$ as tensor categories. 
It often happens that 
a subsystem $\Delta$ of a finite system of non-degenerate braiding 
$\hat{\Delta}$ is degenerate. A typical example is $\hat{\Delta}=$ a full 
system of WZW $SU(N)_k$-model and $\Delta=$ the grading 0 part of 
$\hat{\Delta}$. (The grading is introduced by the cyclic group 
$\mathbb{Z}_N$ acting on the set of integrable highest weight modules of 
level $k$. \cite{Kohno-Takata}) 
In the case of 
$SU(2)_k$ and $SU(3)_k$, they succeeded to describe the 
quantum double of $\Delta$ in terms of $\hat{\Delta}$ \cite{EK2}. Later, by using 
sector theory, Izumi obtained the quantum double of $\Delta$ in 
the case of $SU(N)_k$ which description is quite close to M\"uger's crossed product, 
namely dividing the double category $\hat{\Delta} \otimes \hat{\Delta}^{op}$ by 
the group symmetry $\mathbb{Z}_N$. In this paper, we will generalize Izumi's 
argument to obtain the description of the quantum double of $\Delta$ in the language 
of M\"uger's crossed product, under the assumption that $\hat{\Delta}$  
is a minimal non-degenerate extension i.e.,  $\hat{\Delta} \cap \Delta'=\Delta \cap 
\Delta'$.

M\"uger's theory of crossed product has its origin at a conjecture by 
Rehren \cite{Re1}: Extending endomorphisms on the observable algebra to the 
ones on the field algebra removes the degeneracy of the braiding. M\"uger solved this 
conjecture in \cite{Muger0} and he noticed that it could be possible to formulate the 
whole theory in terms of tensor category \cite{Muger2}.  His formulation 
crucially depends on Doplicher-Roberts duality theory \cite{DR1, DR2}. 
(See \cite{Yamagami} for another equivalent approach to crossed products.)
It should be mentioned that at almost the same time as M\"uger's work, Brugui\`eres 
developed how to construct a modular category out of a certain ribbon 
category,  based on Deligne's internal characterization of tannakien category in 
characteristic 0. 

These modularizations have obvious applications to the 
Reshetikhin-Turaev TQFT. Brugui\`eres himself examined some cases such as 
$SL(N)$ and $PSL(N)$ as examples \cite{Br}. Sawin used M\"uger's 
machinery to obtain a modular category out of  closed subsets of the 
Weyl alcove of a simple Lie algebra, which is essentially the case 
dividing some ribbon categories associated with simple Lie algebras by 
the cyclic group actions. He also obtained a topological invariant of closed 
3-manifolds associated with such modular categories \cite{Sawin}.

Since we have Longo-Rehren inclusions $A \supset B_\Delta \supset 
B_{\hat{\Delta}}$ for a minimal non-degenerate extension $\hat{\Delta} \supset \Delta$, 
we can construct the Reshetikhin-Turaev invariant from the data of the 
quantum double of $\Delta$. As an 
application of an orbifold aspect of the inclusions $A \supset B_\Delta \supset 
B_{\hat{\Delta}}$, we will have a simpler description of the Reshetikhin-Turaev 
invariant of closed 3-manifolds constructed from the quantum double of 
$\Delta$.

This article is organized as follows. \\
In Sec.2, we collect some terminologies we need in this article. In particular, we will 
make quick (and somewhat brutal) reviews on the $\alpha$-induction and M\"uger's 
theory of crossed product. In Sec.3, we explicitly compute the $\alpha$-induction 
for the subfactor $M \rtimes \hat{G} \supset M$, where $M \rtimes 
\hat{G}$ is the crossed product factor by group dual. We will prove that 
the subsectors of $\alpha$-induced sectors for $M \rtimes \hat{G} 
\supset M$ forms a modular category. This result is a folklore 
among the experts, but to the best of author's knowledge, there is no 
description in the present fashion. In Sec.4, we construct the Longo-Rehren inclusions 
$A \supset B_\Delta \supset B_{\hat{\Delta}}$ from a minimal 
non-degenerate extension 
$\hat{\Delta} \supset \Delta$, and we will prove that $B_\Delta
\supset B_{\hat{\Delta}}$ is conjugate to $B_{\hat{\Delta}} \rtimes \hat{G} \supset 
B_{\hat{\Delta}}$. This implies that the quantum double of $\Delta$ can be described by 
$\hat{\Delta}$ and M\"uger's crossed product. In Sec.5, we will apply 
the result obtained in Sec.4 to the Reshetikhin-Turaev invariant 
constructed from the quantum double of $\Delta$. 
Combined with the result in \cite{KSW}, we can have the statement that 
the Turaev-Viro-Ocneanu invariant constructed from $\Delta$ is described 
by the sum of the product of the framed link invariant constructed from 
$\hat{\Delta}$ and its complex conjugate, which gives a special case of Ocneanu's 
theorem \cite{Ocneanu}. 

\medskip
\noindent
{\bf Acknowledgment} The author would like to thank Professors 
Y. Kawahigashi and M. Izumi for valuable discussions in the early stage 
of this work. He also thanks Professor Yamagami and Dr. M\"uger for 
comments on the preliminary version of this manuscript.

\section{Preliminaries}
\subsection{Braided system of endomorphisms}
{\it Basics on sector theory.}\ 
Let $M$, $N$ be infinite factors, and we denote by ${\rm Mor}(N,M)_0$ 
the set of unital normal $*$-homomorphisms from $N$ to $M$ whose image 
has a finite index. The statistical dimension $d(\rho)$ of $\rho \in 
{\rm Mor}(N,M)_0$ is given by $d(\rho)=[M:\rho(N)]^{1/2}$. 
$\rho \in {\rm Mor}(N,M)_0$ is called {\it 
irreducible} if $M \cap \rho(N)' \cong \mathbb{C} 1_M$. For $\rho, 
\sigma \in {\rm Mor}(N,M)_0$, the intertwiner space ${\rm Hom}(\rho, 
\sigma)$ is defined by ${\rm Hom}(\rho, \sigma)=\{ V \in M| V \rho(x)
=\sigma(x)V, \; x \in N \}$. For every $\rho \in {\rm Mor}(N,M)_0$, 
there are isometries $R_\rho \in {\rm Hom}(id, \bar{\rho} \rho)$, 
$\bar{R}_\rho \in {\rm Hom}(id, \rho\bar{\rho})$ satisfying 
$\bar{R}_\rho^* \rho(R_\rho)=\frac{1}{d(\rho)}$, $R_\rho^* \bar{\rho}
(\bar{R}_\rho)=\frac{1}{d(\rho)}$. $\bar{\rho}$ is called the 
{\it conjugate} of $\rho$. 

The unitary equivalence class of $\rho \in {\rm Mor}(N,M)_0$ is called 
a {\it sector}, and we denote by $[\rho]$ the sector of $\rho$. For 
sectors $[\rho_1], [\rho_2], [\rho]$, $\rho_1, \rho_2, \rho \in 
{\rm Mor}(M,M)_0$, we have the product $[\rho_1][\rho_2]=[\rho_1 \cdot 
\rho_2]$, the conjugation $\overline{[\rho]}=[\bar{\rho}]$ and the 
direct sum: $[\rho_1] \oplus [\rho_2]=[\rho]$, where $\rho(x)=t_1^*
\rho_1(x)t_1+t_2^*\rho_2(x)t_2$ and $t_1, t_2 \in M$ satisfying 
$t_i^*t_j=\delta_{i,j}1$ and $t_1t_1^*+t_2t_2^*=1$. 

Let $M \supset N$ be an inclusion of infinite factors with finite 
index $\lambda=[M:N]$ and $\gamma$ be its canonical endomorphism. Then, 
it is known that there exist isometries $v \in {\rm Hom}(id, \gamma)$ and 
${\rm Hom}(\gamma, \gamma^2)$ satisfying 
\begin{equation}\label{vw}
 v^* w=w^* \gamma(v)=\lambda^{-1/2}1, \ w^*\gamma(w)=ww^*, \ 
\gamma(w)w=w^2.
\end{equation}
Moreover, $M=Nv$ and the conditional expectation $E$ from $M$ onto $N$ 
is given by $E(x)=w^*\gamma(x)w$. \\

\noindent
{\it Braided system of endomorphisms.} 
Let $M$ be an infinite factor, and $\Delta_0$ be a system of 
irreducible endomorphisms in ${\rm End}(M)_0={\rm Mor}(M,M)_0$. More 
specifically, $\Delta_0$ is the set of irreducible 
normal $*$-endomorphisms of $M$ closed under the following sector operations:\\
\[
\begin{array}{cl}
{\rm (i)}&  \text{Different elements in }\Delta_0\text{ are inequivalent.} \\
{\rm (ii)}& id_M \in \Delta_0. \\
{\rm (iii)}&  \text{For every }\xi \in \Delta_0\text{ there exists }\bar{\xi} \in 
\Delta_0\text{ such that }\overline{[\xi]}=[\bar{\xi}]. \\
{\rm(iv)} & \text{There exists a non-negative integer }N_{\xi \eta}^\zeta 
\text{ such that }[\xi][\eta]=\oplus_{\zeta \in \Delta_0} N_{\xi \eta}^\zeta [\zeta]. 
\end{array}
\]
\noindent
We denote by $\Delta$ the subset of ${\rm End}(M)_0$ whose element is 
decomposed into finite direct sums of the elements in $\Delta_0$ as sectors. 

A system of endomorphisms $\Delta_0$ is called {\it braided} if for any 
$\lambda, \mu \in \Delta_0$ there exists a unitary intertwiner 
$\varepsilon(\lambda, \mu) \in {\rm Hom}(\lambda \cdot \mu, \mu \cdot \lambda)$ 
with $\varepsilon(id, \mu)=\varepsilon(\lambda, id)=1$ satisfying the following 
(the Braiding-Fusion equations):\\
For any $\lambda, \mu, \nu \in \Delta_0$, $t \in {\rm Hom}(\lambda, \mu \cdot 
\nu)$, 
\begin{eqnarray}
\sigma(t) \varepsilon(\lambda, \sigma)=\varepsilon(\mu, \sigma) \mu(\varepsilon
(\nu, \sigma)) t \label{BFE1}\\
t\varepsilon(\sigma, \lambda)=\mu(\varepsilon(\sigma, \nu)) \varepsilon(\sigma, \mu) 
\sigma(t) \\
\sigma(t)^* \varepsilon(\mu, \sigma) \mu(\varepsilon(\nu, \sigma))=\varepsilon
(\lambda, \sigma)t^* \\
t^* \mu(\varepsilon(\sigma, \nu))\varepsilon(\sigma, \mu)=\varepsilon(\sigma, 
\lambda)\rho(t)^*. \label{BFE4}
\end{eqnarray}
We call above $\varepsilon$ a {\it braiding} on $\Delta_0$. For a given braiding 
$\varepsilon(\lambda, \mu)$ on $\Delta_0$, unitary intertwiners $\varepsilon
(\mu, \lambda)^*$ also satisfies the above conditions of the braiding. 
We will use the notations  $\varepsilon^+(\lambda, \mu)=\varepsilon(\lambda, 
\mu)$ and $\varepsilon^- (\lambda, \mu)=\varepsilon(\mu, \lambda)^*$ to 
emphasize the difference. 

A braiding on $\Delta_0$ can be extended to a braiding on $\Delta$ in 
the following way:\\
\[
\varepsilon(\lambda, \mu)=\sum_{i=1}^n \sum_{j=1}^m s_j \mu_j(t_i) 
\varepsilon(\lambda_i, \mu_j) \lambda(s_j^*)t_i^*,  
\]
where $\lambda(x)=\sum_{i=1}^n t_i \lambda_i(x)t_i^*$ and 
$\mu(x)=\sum_{j=1}^m s_j \mu_j(x)s_j^*$ (i.e., $\lambda=\oplus_{i=1}^n 
\lambda_i$ and $\mu=\oplus_{j=1}^m \mu_j$).

\noindent
{\it Degenerate sectors.} 
A sector $\xi \in \Delta$ is said to be {\it degenerate} if $\varepsilon^+
(\xi, \eta)=\varepsilon^-(\xi, \eta)$ for every $\eta \in \Delta_0$. 
$\Delta$ is said to be {\it non-degenerate} if $id_M$ is the only degenerate 
sector. We denote the set of all of degenerate sectors in $\Delta$ by 
$\Delta^d$ and the set of all of irreducible sectors in $\Delta^d$ by 
$\Delta^d_0$. 
Note that $\Delta^d$ is a symmetric $C^*$-tensor subcategory of 
$\Delta$ with direct sums, subobjects and conjugates. 

For $\xi \in \Delta^d_0$, $\phi_\xi(\varepsilon(\xi,\xi))
=\lambda_\xi \in \mathbb{C}$, where $\phi_\xi$ is the standard 
left inverse of $\xi$. The polar decomposition of $\lambda_\xi$ is 
given by $\frac{\omega_\xi}{d(\xi)}$. It is easy to show that 
$\omega_\xi=\pm 1$ for $\xi \in \Delta^d$ (more generally, for an 
object in a symmetric $C^*$-tensor category). $\Delta^d$ is said to be 
{\it even} if $\omega_\xi=1$ for every irreducible $\xi \in \Delta^d$. 
We assume  $\Delta^d$ is even in the sequel. Then, by Doplicher-Roberts 
duality theory, there exists a finite group $G$ up to isomorphism such 
that $\Delta^d \cong \hat{G}$, where $\hat{G}$ is a category of finite 
dimensional unitary representations of $G$. 

\noindent
$\alpha$-{\it induction.} Let $M \supset N$ be an inclusion of infinite 
factors with finite index 
and $\gamma$ be its canonical endomorphism. Let $\Delta_0 \subset {\rm 
End}(N)_0$ be a braided system of endomorphisms with a braiding 
$\varepsilon$. We define the $\alpha$-induced endomorphism of $\lambda 
\in \Delta_0$ $\alpha_\lambda \in {\rm End}(M)$  by 
\[
 \alpha_\lambda =\gamma^{-1} \cdot Ad(\varepsilon(\lambda, \theta)) 
 \cdot \lambda \cdot \gamma, 
\]
where $\theta=\gamma|_N$.

The systematic use of $\alpha$-induction was first made by Xu \cite{Xu}, 
and further studied in a series of papers by B\"ockenhauer and Evans 
\cite{BE1, BE2, BE3}. 
We list some properties of the $\alpha$-induction \cite{BE1, Xu}:
\[
\begin{array}{cl}
 {\rm (i)}& d(\alpha_\lambda)=d(\lambda)  \\
 {\rm (ii)}& \alpha_\lambda \cdot \alpha_\mu=\alpha_{\lambda \cdot \mu}\  \text{ 
 for any}\   \lambda, \mu \in \Delta_0 \\
{\rm (iii)}& \alpha_\mu \cdot \alpha_\lambda=Ad(\varepsilon(\lambda, \mu)) 
 \cdot \alpha_\lambda \cdot \alpha_\mu \  \text{ 
 for any}\   \lambda, \mu \in \Delta_0 \\
 {\rm (iv)}& {\rm If}\  [\lambda]=[\lambda_1] \oplus [\lambda_2], \ \lambda, 
 \lambda_1, \lambda_2 \in \Delta, \ {\rm then}\ [\alpha_\lambda]=[\alpha_{\lambda_1}]
 \oplus [\alpha_{\lambda_2}] \ {\rm and}  \\
 {\rm (v)}&  [\alpha_{\bar{\lambda}}]=[\overline{\alpha_\lambda}], \lambda 
 \in \Delta_0.
\end{array}
\]
The $\alpha$-induction on $\Delta_0$ is extended to the one on 
$\Delta$ preserving the above properties.

\subsection{Premodular categories and M\"uger's crossed product}

To define M\"uger's crossed product, we need some terminologies 
from category theory. See \cite{LR} for the basics on $C^*$-tensor 
category and \cite{Muger1} for the full description of crossed 
product.

\begin{asmp}\label{assumption}
We assume that $\cal C$ is a $C^*$-tensor category with conjugate, 
direct sums, subobjects, irreducible unit object $\iota$ and a 
unitary braiding $\varepsilon$. 
\end{asmp}
We use the following notations which are popular in the context of 
the algebraic quantum field theory: \\
We use small Greek letters $\rho, \sigma$ etc for objects of $\cal C$, 
and the tensor product is denoted by $\rho \sigma$ 
instead of $\rho \otimes \sigma$. For operations of arrows, we denote 
the composition of arrows $S \in {\rm Hom}(\rho, \sigma)$, $T \in 
{\rm Hom}(\sigma, \tau)$ by $T \circ S \in {\rm Hom}(\rho, \tau)$, 
the tensor product of $S \in {\rm Hom}(\rho_1, \sigma_1)$, $T \in 
{\rm Hom}(\rho_2, \sigma_2)$ by $S \times T \in {\rm Hom}(\rho_1 \rho_2, 
\sigma_1 \sigma_2)$. We denote by ${\cal C}_0$ the set of isomorphism 
classes of irreducible objects.

We remark that under Assumption \ref{assumption} $\cal C$ is 
a ribbon category and we denote a twist for each irreducible object $\rho \in 
{\cal C}$ by $\omega_\rho$. 

Since we assume that $\cal C$ has a conjugate $\bar{\rho}$ for each 
object $\rho$, there are $R_\rho \in {\rm Hom}(\iota, \bar{\rho} \rho)$ 
and $\bar{R}_\rho \in {\rm Hom}(\iota, \rho \bar{\rho})$ satisfying 
\[
 {\bar{R}_\rho}^* \times id_\rho \circ id_\rho \times R_\rho =id_\rho, \ 
{R_\rho}^* \times id_\rho \circ id_{\bar{\rho}} \times \bar{R}_\rho = id_\rho.
\]
Then, the dimension of an irreducible object $\rho$ is defined by $d(\rho)=
{R_\rho}^* \circ R_\rho$, which takes its value in $[1, \infty)$. 
(This definition of the dimension extends to reducible objects.) 

If the set ${\cal C}_0$ is finite, the category is called {\it rational}. 
Then, its dimension is defined by $\dim {\cal C}=\sum_{\xi \in {\cal C}_0} 
d(\xi)^2$. In subfactor context, this is called the global index. 

When ${\cal C}$ is rational, then we set the complex number 
\[
 S'(\xi, \eta) id_\iota= ({R_\xi}^* \times {\bar{R}_\eta}^*) 
\circ (id_{\bar{\xi}} \times \varepsilon(\eta, \xi) \circ 
\varepsilon(\xi, \eta) \times id_{\bar{\eta}}) \circ (R_\xi \times \bar{R}_\eta) 
\]
for $\xi, \eta \in {\cal C}_0$. 
One can prove that $S'(\xi, \eta)$ does not depend on the choice of 
representatives of $\xi$ and $\eta$. 

If $S'$ is invertible, ${\cal C}$ is called {\it modular}. When ${\cal C}$ 
is modular, the matrices 
\[
 S={\dim {\cal C}}^{-\frac{1}{2}} S', \ 
T=\left( \frac{\Delta_{\cal C}}{|\Delta_{\cal C}|} \right)^{\frac{1}{3}}
{\rm Diag}(\omega_\xi) 
\]
are unitaries and satisfy the relations 
\[
 S^2=(ST)^3=C,\ TC=CT,
\]
where $\Delta_{\cal C}=\sum_{\xi \in {\cal C}_0} d(\xi)^2 
\omega(\xi)^{-1}$ and $C=\delta_{\xi, \bar{\eta}}$. 

\begin{defn}
If ${\cal C}$ satisfies Assumption \ref{assumption} and is rational, 
we say ${\cal C}$ is $C^*$-premodular. 
\end{defn}
For a $C^*$-premodular category ${\cal C}$ and its full subcategory 
${\cal S}$, we define ${\cal C} \cap {\cal S}'$, a full subcategory of 
${\cal C}$, by ${\rm Obj}\; {\cal C} \cap {\cal S}' =\{ \rho \in 
{\cal C} | \varepsilon(\sigma, \rho) \circ \varepsilon(\rho, \sigma)=
id_{\rho \sigma} \ {\rm for \ all}\  \sigma \in {\cal S}\}$. We remark 
that if $\cal C$ is modular we have $\dim \; {\cal C} \cap {\cal S}'=
\frac{\dim \; {\cal C}}{\dim \; {\cal S}}$ by Theorem 3.2 \cite{Muger3}. 

Let ${\cal C}$ be a $C^*$-premodular category and we set ${\cal D}_{\cal 
C}={\cal C} \cap {\cal C}'$. We assume that ${\cal D}_{\cal C}$ is even, 
i.e., twist $\omega_{\xi}=1$ for each irreducible object $\xi$. 
Then, by Doplicher-Roberts duality theory \cite{DR1, DR2}, there is 
a finite group such that ${\cal D}_{\cal C}$ is equivalent to $U(G)$ 
as symmetric tensor $*$-categories with conjugates, where $U(G)$ is a 
category of finite dimensional unitary representations of $G$. In the 
following, we use the symbol $\boxtimes$ for the tensor product of 
$U(G)$. 

Let $F$ be an invertible functor from ${\cal D}_{\cal C}$ to $U(G)$ 
which gives the equivalence, $\hat{G}$ be the set of all isomorphism 
classes of irreducible objects in ${\cal D}_{\cal C}$, 
$\{ \gamma_k | k \in \hat{G}\}$ be a section of objects in ${\cal 
D}_{\cal C}$ such that $\gamma_0 = \iota$ and 
${\cal H}_k=F(\gamma_k)$. We choose an orthonormal basis  
$\{ V_{k,l}^{m,\alpha} \}_{\alpha=1}^{N_{kl}^m}$ of 
${\rm Hom}(\gamma_m, \gamma_k \gamma_l)$. 

Then, a category ${\cal C} \rtimes_0 {\cal D}_{\cal C}$ is defined 
in the following manner.
\begin{itemize}
\item ${\rm Obj}\; {\cal C} \rtimes_0 {\cal D}_{\cal C}= 
{\rm Obj}\; {\cal C} $ with the same tensor product as 
${\cal C}$ 
\item ${\rm Hom}_{{\cal C} \rtimes_0 {\cal D}_{\cal C}} (\rho, \sigma)
=\bigoplus_{k \in \hat{G}} {\rm Hom}_{{\cal C}} (\gamma_k\rho, \sigma) 
\otimes {\cal H}_k$.
\end{itemize}
Let $k, l \in \hat{G}$, $S \otimes \psi_k \in {\rm Hom}_{{\cal C} 
\rtimes_0 {\cal D_{\cal C}}}(\rho, \sigma)$ and $T \otimes \psi_l \in 
{\rm Hom}_{{\cal C} \rtimes_0 {\cal D_{\cal C}}}(\sigma, \tau)$, where 
$T \in {\rm Hom}(\gamma_l \rho, \sigma)$, 
$S \in {\rm Hom}(\gamma_k \sigma, \tau)$ and $\psi_k \in {\cal H}_k$, 
$\psi_l \in {\cal H}_l$. We define the composition of arrows 
$S \otimes \psi_k \circ T \otimes \psi_l \in {\rm Hom}_{{\cal C} 
\rtimes_0 {\cal D_{\cal C}}}(\rho, \tau)$ by 
\begin{equation}\label{cp1}
 S \otimes \psi_k \circ T \otimes \psi_l 
=\bigoplus_{k \in \hat{G}} \sum_{\alpha=1}^{N_{kl}^m} 
S \circ id_{\gamma_k} \times T \circ V_{k,l}^{m,\alpha} \times id_\rho 
\otimes F(V_{k,l}^{m,\alpha})^* (\psi_k \boxtimes \psi_l)
\end{equation}
and extend this linearly.

Let $k, l \in \hat{G}$, $S \otimes \psi_k \in {\rm Hom}_{{\cal C} 
\rtimes_0 {\cal D_{\cal C}}}(\rho_1, \sigma_1)$ and $T \otimes \psi_l \in 
{\rm Hom}_{{\cal C} \rtimes_0 {\cal D_{\cal C}}}(\rho_2, \sigma_2)$, 
where $S \in {\rm Hom}(\gamma_k \rho_1, \sigma_1)$, 
$T \in {\rm Hom}(\gamma_l \rho_2, \sigma_2)$ and $\psi_k \in {\cal H}_k$, 
$\psi_l \in {\cal H}_l$. We define the tensor product of arrows 
$S \otimes \psi_k \times T \otimes \psi_l \in {\rm Hom}_{{\cal C} 
\rtimes_0 {\cal D_{\cal C}}}(\rho_1 \rho_2, \sigma_1 \sigma_2)$ by
\begin{eqnarray}\label{cp2}
S \otimes \psi_k \times T \otimes \psi_l 
= &&\bigoplus_{k \in \hat{G}} \sum_{\alpha=1}^{N_{kl}^m} 
S \times T \circ id_{\gamma_k} \times \varepsilon(\gamma_l, \rho_1) 
\times id_{\rho_2} \circ V_{k,l}^{m,\alpha} \times id_{\rho_1 \rho_2} 
\nonumber \\
&\otimes& F(V_{k,l}^{m,\alpha})^*(\psi_k \boxtimes \psi_l)
\end{eqnarray}
and extend this linearly.

Let $S \otimes \psi \in {\rm Hom}_{{\cal C} \rtimes_0 {\cal D_{\cal C}}}
(\rho, \sigma)$, where $S \in {\rm Hom}(\gamma_k \rho, \sigma)$ and 
$\psi_k \in {\cal H}_k$. 
We define the $*$-operation of the arrows $(S \otimes \psi)^* \in 
{\rm Hom}_{{\cal C} \rtimes_0 {\cal D_{\cal C}}}(\sigma, \rho)$ by 
\begin{equation}\label{cp3}
(S \otimes \psi)^*={R_k}^* \times id_\rho \circ id_{\bar{\gamma_k}} 
\times S^* \otimes \langle \psi_k \boxtimes \cdot, F(\bar{R}_k) \Omega \rangle,
\end{equation}
where $\Omega$ is a unit vector in the trivial representation ${\cal H}_0 \cong 
\mathbb{C}$ such that $\Omega \boxtimes \psi = \psi \boxtimes 
\Omega=\psi$ for all $\psi \in {\rm Obj}\; U(G)$. 

It turns out that ${\cal C} \rtimes_0 {\cal D}_{\cal C}$ is a 
$C^*$-tensor category with conjugates and direct sums. 

\begin{rem}
 For ${\cal C}$, we have another braiding $\varepsilon^-(\lambda, \mu)=
\varepsilon(\mu,\lambda)^*$. When we need 
 to clarify which braiding we used, we will write ${\cal C} \rtimes_{0, +} 
{\cal D}_{\cal C}$ and ${\cal C} \rtimes_{0, -} {\cal D}_{\cal C}$ 
depending on the choice of the braiding $\varepsilon^+$ and 
$\varepsilon^-$, respectively. 
\end{rem}

${\cal C} \rtimes_0 {\cal D}_{\cal C}$ is not closed under subobjects in general. 
However, we can enlarge ${\cal C} \rtimes_0 {\cal D}_{\cal C}$ to be closed under 
the subobjects. Such a procedure is called {\it closure} in Definition 
3.11 in \cite{Muger2}. We denote the closure of ${\cal C} \rtimes_0 
{\cal D}_{\cal C}$ by ${\cal C} \rtimes {\cal D}_{\cal C}$ and call it 
the {\it crossed product} of ${\cal C}$ by ${\cal D}_{\cal C}$.  We remark 
that $\dim \; {\cal C} \rtimes {\cal D}_{\cal C}=\frac{\dim \; {\cal C}}{\dim \; 
{\cal D}_{\cal C}}$. 

It is important to mention that ${\cal C} \rtimes {\cal D}_{\cal C}$ is 
a modular category due to Theorem 4.4 in \cite{Muger2}. 

\begin{rem}
 M\"uger has constructed a $C^*$-tensor category 
${\cal C} \rtimes {\cal S}$ with conjugates , direct sums and subobjects 
from $\cal C$ (not rational in general) and $\cal S$, where $\cal S$ is a 
symmetric $C^*$-tensor subcategory of $\cal C$, not necessarily ${\cal S} \subset 
{\cal D}_{\cal C}$. See \cite{Muger2} for details. 
\end{rem}
\section{M\"uger's crossed product and $\alpha$-induction for subfactors}

Let $M$, $\Delta$  and $\Delta^d$ be as in Subsection 2.1, and we 
assume that $\Delta_0$ is a finite set. We further assume that 
$\Delta^d$ is even and $\Delta^d \cong U(G)$, where $G$ is a finite 
group. Then, by Doplicher-Roberts duality theory there exists 
 a factor, denoted by $M \rtimes \hat{G}$, 
which contains $M$ as a subfactor with index $|G|$. See 
\cite{DR1} and \cite{DR2} for the detailed accounts.

We may assume that $M \rtimes \hat{G}$ is generated by $M$ and 
isometries $\{\psi^{(\sigma)}_i$, $i=1, \cdots, d(\sigma), \sigma \in \
\Delta^d_0\}$ satisfying : 
 \begin{eqnarray}
& \psi^{(\iota)}:=\psi_1^{(\iota)}=1  \\
&  {\psi^{(\sigma)}_i}^* \psi^{(\sigma')}_j=\delta_{i,j}\delta_{\sigma, 
  \sigma'} \label{r1} \\
&  \sum_{i=1}^{d(\sigma)} \psi^{(\sigma)}_i 
 {\psi^{(\sigma)}_i}^*=1 \label{r2} \\
& \psi^{(\sigma)}_i x=\sigma(x) \psi^{(\sigma)}_i, 
\ x \in M \label{r3} \\
& \psi^{(\rho)}_i \psi^{(\sigma)}_j=\sum_{\tau \in \Delta^d_0} 
\sum_{k=1}^{d(\tau)} V_{(\rho,i) (\sigma,j)}^{(\tau, k)} 
\psi^{(\tau)}_k  \label{r4} \\
& {\psi^{(\sigma)}_i}^*=R_\sigma^* \psi^{(\bar{\sigma})}_i \label{r5} \\
& \sum_{i=1}^{d(\sigma_1)} \sum_{j=1}^{d(\sigma_2)} \psi^{(\sigma_2)}_j 
\psi^{(\sigma_1)}_i {\psi^{(\sigma_2)}_j}^* {\psi^{(\sigma_1)}_i}^*=\varepsilon
(\sigma_1, \sigma_2), \label{r6}
 \end{eqnarray}
where $V_{(\rho,i) (\sigma,j)}^{(\tau, k)} \in {\rm Hom}(\tau, 
\rho \cdot \sigma)$ and $R_\sigma \in {\rm Hom}(\iota, \bar{\sigma} 
\cdot \sigma)$.
\begin{rem}
{\rm (1)} It is known that $\{\psi^{(\sigma)}_i$, $i=1, \cdots, d(\sigma), \sigma 
\in \ \Delta^d_0\}$ is a left $M$-module basis. \\
{\rm (2)} When $x=\sum_{\sigma,i}t_i^{(\sigma)} \psi_i^{(\sigma)} \in M 
\rtimes \hat{G}$, the conditional expectation $E:M \rtimes \hat{G} \longrightarrow M$ 
is given by $E(x)=t^{(\iota)}$. By computations, one has $E(\psi_i^{(\sigma)} 
{\psi_j^{(\rho)}}^*)=\delta_{\sigma, \rho}\delta_{i,j}\frac{1}{d(\sigma)}$, 
where $\lambda=[M \rtimes \hat{G}:M]$.
\end{rem}
\begin{lem}\label{v}
Let $v=\sum_{\sigma, i} t_i^{(\sigma)} \psi_i^{(\sigma)} \in 
{\rm Hom}(id, \gamma)$. 
Then, we have the relations $t_i^{(\sigma)}=d(\sigma) E(v 
{\psi_i^{(\sigma)}}^*) \in {\rm Hom}(\sigma, \theta)$ 
and $\psi_i^{(\sigma)}=\frac{\lambda}{d(\sigma)} {t_i^{(\sigma)}}^* v$. 
Furthermore, 
$t_i^{(\sigma)}$, $i=1, \cdots, d(\sigma)$ satisfy ${t_i^{(\sigma)}}^* 
t_j^{(\rho)}=\delta_{\sigma, \rho} \delta_{i,j} 
\frac{d(\sigma)}{\lambda}$ and $\sum_{\sigma, i}\frac{\lambda}{d(\sigma)}
t_i^{(\sigma)} {t_i^{(\sigma)}}^*=1$. 
\end{lem}

\noindent
{\it Proof.}
Applying the conditional expectation $E$ to the equation 
$v {\psi_j^{(\rho)}}^*=\sum_{\sigma, i} t_i^{(\sigma)} \psi_i^{(\sigma)} 
{\psi_j^{(\rho)}}^*$, we have 
\[
 E(v {\psi_j^{(\rho)}}^*)=\sum_{\sigma, i}t_i^{(\sigma)} 
E(\psi_i^{(\sigma)} {\psi_j^{(\rho)}}^*)
= t_j^{(\rho)} \frac{1}{d(\rho)}.
\]
Therefore, $t_i^{(\sigma)}=d(\sigma) E(v {\psi_i^{(\sigma)}}^*)$. 
Multiplying $v$ from the left of the equality ${\psi_i^{(\sigma)}}^* 
\sigma(x)=x {\psi_i^{(\sigma)}}^*$, $x \in M$, we have 
\[
 v {\psi_i^{(\sigma)}}^* \sigma(x) =vx{\psi_i^{(\sigma)}}^* 
=\gamma(x) v {\psi_i^{(\sigma)}}^* =\theta(x) v {\psi_i^{(\sigma)}}^*.
\]
Apply the conditional expectation $E$ to the above equality, then we 
have 
\[
 E(v {\psi_i^{(\sigma)}}^*) \sigma(x)=\theta(x) E(v {\psi_i^{(\sigma)}}^*).
\]
Hence, $t_i^{(\sigma)}=d(\sigma) E(v {\psi_i^{(\sigma)}}^*) \in 
{\rm Hom}(\sigma, \theta)$. 

Let us compute ${t_i^{(\sigma)}}^* t_j^{(\rho)}$. 
\begin{eqnarray*}
 {t_i^{(\sigma)}}^* t_j^{(\rho)} 
&=& d(\sigma) d(\rho) E(\psi_i^{(\sigma)} v^*) E(v {\psi_j^{(\rho)}}^*) \\
&=& d(\sigma) d(\rho) w^* \gamma(\psi_i^{(\sigma)} v^*) ww^* 
\gamma(v {\psi_j^{(\rho)}}^*) w \\
&=& d(\sigma) d(\rho) w^* \gamma(\psi_i^{(\sigma)})\gamma(v^*) w w^* 
\gamma(v)\gamma({\psi_j^{(\rho)}}^*) w  \\
&=& \lambda^{-1} d(\sigma) d(\rho) w^* \gamma(\psi_i^{(\sigma)}{\psi_j^{(\rho)}}^*)w  \\
&=& \lambda^{-1} d(\sigma) d(\rho) E(\psi_i^{(\sigma)}{\psi_j^{(\rho)}}^*) \\
&=& \delta_{\sigma,\rho} \delta_{i,j} \frac{d(\sigma)}{\lambda}.
\end{eqnarray*}

Next, we compute ${t_i^{(\sigma)}}^*v$. 
\begin{eqnarray*}
 {t_i^{(\sigma)}}^* v 
&=& d(\sigma) E(\psi_i^{(\sigma)}v^*)v \\
&=& d(\sigma) w^* \gamma(\psi_i^{\sigma}) \gamma(v^*) w v \\
&=& \lambda^{-1/2} d(\sigma) w^* \gamma(\psi_i^{(\sigma)}) v \\
&=& \lambda^{-1/2} d(\sigma) w^* v \psi_i^{(\sigma)} \\
&=& \frac{d(\sigma)}{\lambda} \psi_i^{(\sigma)}.
\end{eqnarray*}
Hence, $\psi_i^{(\sigma)}=\frac{\lambda}{d(\sigma)}{t_i^{(\sigma)}}^* v$. 

Then, we have the following identity
\[
 v=\sum_{\sigma, i} t_i^{(\sigma)} \psi_i^{(\sigma)}
=\sum_{\sigma, i} \frac{\lambda}{d(\sigma)} t_i^{(\sigma)} 
{t_i^{(\sigma)}}^* v.
\]
Multiplying $v^*$ from the right of the above equality and applying $E$, then we have 
\[
 \sum_{\sigma, i} \frac{\lambda}{d(\sigma)} t_i^{(\sigma)} {t_i^{(\sigma)}}^*
=1.
\]
This completes the proof. \hfill $\Box$

\begin{prop}\label{chi-loc}
 The equation (\ref{r6}) is equivalent to the identity $\varepsilon(\theta, 
 \theta)v^2=v^2$.
\end{prop}

\noindent
{\it Proof.}
By Lemma \ref{v}, $\sqrt{\frac{\lambda}{d(\sigma)}}t_i^{(\sigma)}$ is 
isometry in ${\rm Hom}(\sigma, \theta)$. Hence, $\varepsilon(\theta,\theta)$ is 
given by 
\begin{eqnarray*}
 \varepsilon(\theta, \theta)
&=& \sum_{\sigma, \sigma'} \sum_{i,j} \sqrt{\frac{\lambda}{d(\sigma')}}
t_j^{(\sigma')} \sigma'(\sqrt{\frac{\lambda}{d(\sigma)}}t_i^{(\sigma)}) 
\varepsilon(\sigma, \sigma') \sigma(\sqrt{\frac{\lambda}{d(\sigma')}} 
{t_j^{(\sigma')}}^*) \sqrt{\frac{\lambda}{d(\sigma)}} {t_i^{(\sigma)}}^* \\
&=& \sum_{\sigma,\sigma'} \sum_{i,j} \lambda^2 d(\sigma)^{-1} 
d(\sigma')^{-1} t_j^{(\sigma')} \sigma'(t_i^{(\sigma)}) 
\varepsilon(\sigma, \sigma') \sigma({t_j^{(\sigma')}}^*) {t_i^{(\sigma)}}^*.
\end{eqnarray*} 

Then, we have
\begin{eqnarray}
 \varepsilon(\theta, \theta)v^2
&=& \sum_{\sigma,\sigma',\tau} \sum_{i,j,k}\lambda d(\sigma')^{-1}
t_j^{(\sigma')} \sigma'(t_i^{(\sigma)})\varepsilon(\sigma,\sigma')
\sigma({t_j^{(\sigma')}}^*) \psi_i^{(\sigma)} t_k^{(\tau)} \psi_k^{(\tau)} 
\nonumber \\
&=& \sum_{\sigma,\sigma',\tau} \sum_{i,j,k}\lambda d(\sigma')^{-1}
t_j^{(\sigma')} \sigma'(t_i^{(\sigma)})\varepsilon(\sigma,\sigma')
\sigma({t_j^{(\sigma')}}^* t_k^{(\tau)}) \psi_i^{(\sigma)}  \psi_k^{(\tau)}
\nonumber \\
&=& \sum_{\sigma,\sigma'} \sum_{i,j} t_j^{(\sigma')} \sigma'(t_i^{(\sigma)})
\varepsilon(\sigma,\sigma') \psi_i^{(\sigma)}\psi_j^{(\sigma')}. \label{evv}
\end{eqnarray}
The equation (\ref{r6}) is equivalent to 
$\varepsilon(\sigma,\sigma')\psi_i^{(\sigma)} \psi_j^{(\sigma')}=\psi_j^{(\sigma')}
\psi_i^{(\sigma)}$, and with this, the formula (\ref{evv}) is equal to
\[
\sum_{\sigma,\sigma'} \sum_{i,j} t_j^{(\sigma')} \sigma'(t_i^{(\sigma)})
 \psi_j^{(\sigma')}\psi_i^{(\sigma)}\\
= \sum_{\sigma,\sigma'} \sum_{i,j} t_j^{(\sigma')}\psi_j^{(\sigma')} t_i^{(\sigma)}
\psi_i^{(\sigma)} \\
= v^2. 
\]

On the contrary, assume that $\varepsilon(\theta, \theta)v^2=v^2$. 
Multiplying $\sigma'(t_i^{(\sigma)})^* {t_j^{(\sigma')}}^*$ from the 
left of $\varepsilon(\theta, \theta)v^2=\sum_{\sigma, \sigma'} 
\sum_{i,j} t_j^{(\sigma')} \sigma'(t_i^{(\sigma)})
\varepsilon(\sigma,\sigma') \psi_i^{(\sigma)}\psi_j^{(\sigma')}$, 
we have 
\[
\sigma'(t_i^{(\sigma)})^* {t_j^{(\sigma')}}^* \varepsilon(\theta, \theta)v^2
=\frac{d(\sigma')d(\sigma)}{\lambda^2} \varepsilon(\sigma,\sigma')
\psi_i^{(\sigma)}\psi_j^{(\sigma')}.
\]
On the other hand, multiplying $\sigma'(t_i^{(\sigma)})^* {t_j^{(\sigma')}}^*$ from the 
left of $v^2$, we have 
\[
 \sigma'(t_i^{(\sigma)})^* {t_j^{(\sigma')}}^* v^2
=\frac{d(\sigma')d(\sigma)}{\lambda^2} \psi_j^{(\sigma')}\psi_i^{(\sigma)}.
\]
Thus, $\varepsilon(\sigma,\sigma')\psi_i^{(\sigma)}\psi_j^{(\sigma')}
=\psi_j^{(\sigma')}\psi_i^{(\sigma)}$. \hfill $\Box$

\begin{rem}
The identity $\varepsilon(\theta,\theta)v^2=v^2$ is called the {\rm chiral 
locality} condition in \cite{BE3}. 
\end{rem}
\begin{lem}\label{alpha-ind}
 For $\lambda \in \Delta$, we have 
\begin{equation}\label{alphaind} \alpha^\pm_\lambda(\psi^{(\sigma)}_i)
=\varepsilon^\pm(\lambda,\sigma)^* \psi^{(\sigma)}_i, 
\end{equation} 
where $\sigma \in {\Delta^d}_0, i=1,  \cdots, d(\sigma)$. 
In particular, $\alpha_\lambda^+=\alpha_\lambda^-$ 
for $\lambda \in \Delta \cap {\Delta^d}'=\{\rho_\xi \in \Delta | 
\varepsilon(\xi, \sigma) \varepsilon(\sigma, \xi)=1, \; \forall \sigma 
\in \Delta^d_0 \}$. 
\end{lem}

\noindent
{\it Proof.} 
Let $\gamma$ be the canonical endomorphism of $M \rtimes \hat{G} 
\supset M$ and $\theta$ the restriction of $\gamma$ to $M$. 

Applying $\gamma$ to (\ref{r3}), we have $\gamma(\psi^{(\sigma)}_i)\theta(x)=
\theta \cdot \sigma(x) \gamma(\psi^{(\sigma)}_i)$. Thus, 
$\gamma({\psi^{(\sigma)}_i}^*) \in {\rm Hom}(\theta \cdot \sigma, \theta)$. 
By the Braiding-Fusion equation (\ref{BFE4}), 
\[
\varepsilon^\pm(\lambda,\theta)^* \theta(\varepsilon^\pm(\lambda,\sigma)^*)
\gamma(\psi^{(\sigma)}_i)
=\lambda(\gamma(\psi^{(\sigma)}_i))\varepsilon^\pm(\lambda,\theta)^*.
\]
Applying $\gamma^{-1}$, 
\[
 \varepsilon^\pm(\lambda,\sigma)^*\psi^{(\sigma)}_i=\gamma^{-1} \cdot 
Ad(\varepsilon^\pm(\lambda,\theta))\lambda \cdot \gamma(\psi^{(\sigma)}_i)
=\alpha^\pm_\lambda(\psi^{(\sigma)}_i).
\]
The last claim is clear because $\varepsilon^+(\lambda, \sigma)=\varepsilon^-
(\lambda, \sigma)$ for $\lambda \in \Delta \cap {\Delta^d}'$.
\hfill $\Box$

\begin{lem}
 For $\lambda, \mu \in \Delta$, 
\[
 {\rm Hom}(\alpha_\lambda, \alpha_\mu)=\{ \sum_{\sigma \in \Delta^d_0} 
 \sum_{i=1}^{d(\sigma)} t^{(\sigma)}_i \psi^{(\sigma)}_i ; 
 t^{(\sigma)}_i \in {\rm Hom}(\sigma \cdot \lambda, \mu), \ i=1, \cdots, d(\sigma), \sigma \in \Delta^d_0 \}.
\]
\end{lem}

\noindent
{\it Proof.}
Let $t \in {\rm Hom}(\alpha_\lambda, \alpha_\mu)$. We may write $t=
\sum_{\sigma \in \Delta^d_0} \sum_{i=1}^{d(\sigma)} t^{(\sigma)}_i \psi^{(\sigma)}_i$.
We remark that this expression is unique. 
For $x \in M$, we have
\[
 \sum_{\sigma, i} t^{(\sigma)}_i \psi^{(\sigma)}_i 
 \lambda(x)=\sum_{\sigma,i} \mu(x) t^{(\sigma)}_i\psi^{(\sigma)}_i.
\]
Since $\psi^{(\sigma)}_i \in {\rm Hom}(id, \sigma)$, the 
above equality is 
\[
 \sum_{\sigma, i} t^{(\sigma)}_i \sigma \cdot \lambda(x) \psi^{(\sigma)}_i 
 =\sum_{\sigma,i} \mu(x) t^{(\sigma)}_i\psi^{(\sigma)}_i.
\]
Thus, $t^{(\sigma)}_i \sigma \cdot \lambda(x)=\mu(x) t^{(\sigma)}_i$ for 
any $x \in M$,  $i=1, \cdots, d(\sigma)$ and $\sigma \in \Delta^d_0$. 

For $t$ above, let us show $t \alpha_\lambda(\psi^{(\sigma')})
=\alpha_\mu(\psi^{(\sigma')})t$, where $\psi^{(\sigma')}$ is an isometry 
in $\{\psi^{(\sigma)}_i, i=1, \cdots, d(\sigma), \sigma \in \Delta^d_0 \}$.
We will show that the left hand side is equal to the right hand side.

\begin{eqnarray*}
 \sum_{\sigma, i} t^{(\sigma)}_i \psi^{(\sigma)}_i 
\varepsilon(\lambda,\sigma')^* \psi^{(\sigma')} &=& \sum_{\sigma, i} t^{(\sigma)}_i
\sigma(\varepsilon(\lambda,\sigma')^*)\psi^{(\sigma)}_i \psi^{(\sigma')} \\
&=& \sum_{\sigma, i} \varepsilon(\mu,\sigma')^* \sigma'(t^{(\sigma)}_i)\varepsilon(\sigma,\sigma')
\psi^{(\sigma)}_i \psi^{(\sigma')} \\
&=& \sum_{\sigma, i} \varepsilon(\mu,\sigma')^* \sigma'(t^{(\sigma)}_i) 
\psi^{(\sigma')} \psi^{(\sigma)}_i \\
&=& \sum_{\sigma, i} \varepsilon(\mu,\sigma')^* \psi^{(\sigma')} t^{(\sigma)}_i \psi^{(\sigma)}_i, 
\end{eqnarray*}
where we used the Braiding-Fusion equation (\ref{BFE1}) for the second equality.
This completes the proof. \hfill $\Box$

\begin{rem}\label{triv-deg}
By the above lemma, we have
\[
 {\rm Hom}(id, \alpha_\rho)=\{  \sum_{i=1}^{d(\rho)} t^{(\rho)}_i \psi^{(\rho)}_i ; 
 t^{(\rho)}_i \rho(x)=\rho(x) t^{(\rho)}_i, \ 
 \forall x \in M, i=1, \cdots, d(\rho) \} 
\] 
for $\rho \in \Delta^d_0$, which is a Hilbert space with 
dimension $d(\rho)$. Since $d(\alpha_\rho)=d(\rho)$, we conclude that 
$\alpha_\rho \cong \oplus_{i=1}^{d(\rho)} id$. This can be read that 
$\alpha$-induction trivializes degenerate sectors. 
\end{rem}

Let $\lambda \in \Delta \cap {\Delta^d}'$ and we use the notation 
$\alpha_\lambda$ instead of $\alpha_\lambda^+=\alpha_\lambda^-$. We denote
by $(\Delta \cap {\Delta^d}')^\alpha$ the subset of ${\rm End}(M \rtimes 
\hat{G})_0$ consisting of subsectors of $\alpha_\lambda$, when $\lambda$ 
varies in $\Delta \cap {\Delta^d}'$.

Thanks to Proposition 
\ref{chi-loc}, we can make a full use of the arguments in the subsection 
3.3 in \cite{BE3}. 
For this, let $\beta, \delta$ be subsectors of $\alpha_\lambda$ and 
$\alpha_\mu$ for some $\lambda, \mu \in \Delta \cap {\Delta^d}'$, respectively. 
We set 
\[
 \varepsilon_r(\beta, \delta)=s^* \alpha_\mu(t^*) \varepsilon(\lambda, \mu)
\alpha_\lambda(s)t \in {\rm Hom}(\beta \cdot \delta, \delta \cdot \beta )
\]
with isometries $t \in {\rm Hom}(\beta, \alpha_\lambda)$, $s \in {\rm Hom}(\delta, 
\alpha_\mu)$. It is proved in Lemma 3.11 \cite{BE3} that 
$\varepsilon_r(\beta, \delta)$ does not depend on $\lambda, \mu$ and on the 
isometries $s, t$. Moreover, $\varepsilon_r(\beta, \delta)$ for 
$\beta, \delta \in (\Delta \cap 
{\Delta^d}')^\alpha$ defines a braiding (called a {\it relative braiding}) on 
$(\Delta \cap {\Delta^d}')^\alpha$ (Corollary 3.13 \cite{BE3}). 

Under these preliminaries, we have the following 
\begin{prop}
 $(\Delta \cap {\Delta^d}')^\alpha$ is a modular category.
\end{prop}

\noindent
{\it Proof.} 
Let $\alpha_\lambda = \oplus_{i=1}^p \beta_i$, $\lambda \in \Delta \cap 
{\Delta^d}'$, and $\delta_j \in (\Delta \cap {\Delta^d}')^\alpha$ such 
that $\alpha_\mu=\oplus_{j=1}^q \delta_j$ for some $\mu \in 
(\Delta \cap {\Delta^d}')^\alpha$. 
Assume $\varepsilon_r(\beta_i, \delta_j)\varepsilon_r(\delta_j, \beta_i)
=1$ for all $j=1, \cdots, q$. Then, we have $\varepsilon(\lambda, \delta)
\varepsilon(\delta, \lambda)=1$ by Lemma 3.14 \cite{BE3}. 
Hence, for $\forall \delta \in \Delta \cap {\Delta^d}'$, we have 
$\varepsilon(\lambda, \delta)\varepsilon(\delta, \lambda)=1$, which 
implies $\lambda \in \Delta^d$. 

Since $\alpha_\lambda = \oplus_{i=1}^{d(\lambda)} id$ by Remark 
\ref{triv-deg}, we have $\beta_i=id$ for all $i=1, \cdots, p$, which  
proves that $\varepsilon_r$ is a non-degenerate braiding on 
$(\Delta \cap {\Delta^d}')^\alpha$. Thus, $(\Delta \cap 
{\Delta^d}')^\alpha$ is modular. \hfill $\Box$

So far, we have discussed the similarities to M\"uger's theory of 
crossed product. In fact, we have the following
\begin{prop}\label{muger}
For the inclusion $M \rtimes \hat{G} \supset M$, 
the image of $\Delta$ by the $\alpha^\pm$-induction is given by $\Delta 
\rtimes_{0,\pm} \Delta^d$. In particular, $(\Delta \cap 
{\Delta^d}')^\alpha$ is naturally identified with 
$(\Delta \cap {\Delta^d}') \rtimes \Delta^d$. 
\end{prop}

\noindent
{\it Proof.}  For the composition of the intertwiners, let 
$s=\sum_{\sigma \in \Delta^d_0} \sum_{i=1}^{d(\sigma)} s^{(\sigma)}_i 
\psi^{(\sigma)}_i \in {\rm Hom}(\alpha_\lambda, \alpha_\mu)$, 
$t=\sum_{\rho \in \Delta^d_0} \sum_{i=1}^{d(\rho)} t^{(\rho)}_i \psi^{(\rho)}_i \in 
{\rm Hom}(\alpha_\mu, \alpha_\nu)$. Then, $ts \in {\rm Hom}(\alpha_\lambda,\alpha_\nu)$ 
defines the composition of morphisms $t$ and $s$. 

We use the notations $\rho=\gamma_k$, $\sigma=\gamma_l$, 
$s^{(l)}=s^{(\gamma_l)}$, $t^{(k)}=t^{(\gamma_k)}$, 
$\psi^{(l)}=\psi^{(\gamma_l)}$ and $\psi^{(k)}=\psi^{(\gamma_k)}$, for 
simplicity. 

It is enough to check the condition for $t=t^{(k)} \psi^{(k)}$ and 
$s=s^{(l)} \psi^{(l)}$ because of linearity.
\begin{eqnarray*}
 t^{(k)} \psi^{(k)} s^{(l)} \psi^{(l)} &=& 
t^{(k)} \gamma_k(s^{(l)}) \psi^{(k)} \psi^{(l)} \\
&=& \sum_{\gamma_m \in \Delta^d_0} \sum_{\alpha=1}^{d(\gamma_m)} 
t^{(k)} \gamma_k(s^{(l)}) V_{k,l}^{m,\alpha} \psi^{(m)}_\alpha, 
\end{eqnarray*}
which is (\ref{cp1}). 

For the tensor product of the intertwiners, let $s \in 
{\rm Hom}(\alpha_{\lambda_1}, \alpha_{\mu_1})$, $t \in {\rm Hom}(\alpha_{\lambda_2}, 
\alpha_{\mu_2})$. Then, $s\alpha_{\lambda_1}(t) \in {\rm Hom}(\alpha_{\lambda_1} 
\cdot \alpha_{\lambda_2}, \alpha_{\mu_1} \cdot \alpha_{\mu_2})$. We 
compute $s\alpha_{\lambda_1}(t)$ in the case 
$s=s^{(k)}\psi^{(k)}$ and $t=t^{(l)}\psi^{(l)}$. 

\begin{eqnarray*}
 s\alpha_{\lambda_1}(t) &=& s^{(k)}\psi^{(k)}t^{(l)}\varepsilon
(\lambda_1,\gamma_l)^* \psi{(l)} \\
&=& s^{(k)}\gamma_k(t^{(l)}) \gamma_k(\varepsilon
(\gamma_{\lambda_1}, \gamma_l)^*) \psi^{(k)} \psi^{(l)} \\
&=& \sum_{\gamma_m \in \Delta^d_0} \sum_{\alpha=1}^{d(\gamma_m)} s^{(k)} 
\gamma_k(t^{(l)}) \gamma_k(\varepsilon(\lambda_1, 
\gamma_l)^*)V_{k,l}^{m,\alpha} \psi^{(m)}_\alpha, 
\end{eqnarray*}  
which is (\ref{cp2}). 

For the $*$-operation of the intertwiners, let 
$t=\sum_{\sigma,i}t^{(\sigma)}_i \psi^{(\sigma)}_i$. We check the 
condition for $t=t^{(\sigma)}\psi^{(\sigma)}$. 
\[
 (t^{(\sigma)}{\psi^{(\sigma)}})^*={\psi^{(\sigma)}}^* {t^{(\sigma)}}^* 
=R_\sigma^* \psi^{(\bar{\sigma})} {t^{(\sigma)}}^* =R_\sigma^* \bar{\sigma}
({t^{(\sigma)}}^*) \psi^{(\bar{\sigma})}, 
\]
which is (\ref{cp3}). 

Thus, the image of the $\alpha$-induction is the crossed product $\Delta 
\rtimes_0 \Delta^d$ in the sense of M\"uger. 

The last claim is immediate from the definitions of 
$(\Delta \cap {\Delta^d}') \rtimes \Delta^d$ and $(\Delta \cap 
{\Delta^d}')^\alpha$. \hfill $\Box$
\section{Longo-Rehren inclusions $A \supset B_\Delta \supset B_{\hat{\Delta}}$}

Let $\Delta$ be a subset of ${\rm End}(M)_0$ with a finite braided 
system $\Delta_0$, $\hat{\Delta} \supset \Delta$ its non-degenerate 
extension. The following definition was first introduced by Ocneanu 
\cite{Ocneanu}.

\begin{defn}
 The non-degenerate extension $\hat{\Delta} \supset \Delta$ is called 
minimal if $\hat{\Delta} \cap \Delta' =\Delta^d$. 
\end{defn}
\noindent
Remark that we have $\dim \hat{\Delta}=\dim \Delta \dim 
\Delta^d$ if the extension is minimal. 

We assume the minimality of the non-degenerate extension $\hat{\Delta} \supset \Delta$ 
in the sequel. 

Let $\{T(^\zeta_{\xi,\eta})_i\}_{i=1}^{N_{\xi,\eta}^\zeta}$ be an 
orthonormal basis of ${\rm Hom}(\zeta, \xi \cdot \eta)$, $\xi,\eta, 
\zeta \in \Delta_0$. Let $M$ be the opposite algebra of $M$ and 
$j:M \longrightarrow M^{op}$ the anti-linear isomorphism. 
We set $A=M \otimes M^{op}$, $\xi^{op}=j \cdot \xi \cdot j$, and 
$\hat{\xi}=\xi \otimes \xi^{op}$. For the isometries $\{V_\xi\}_{\xi \in 
\Delta_0} \subset A$ satisfying $\sum_{\xi \in \Delta_0} V_\xi 
V_\xi^*=1$, we define 
\[
 \gamma_\Delta(x)=\sum_{\xi \in \Delta_0} V_\xi \hat{\xi}(x)V_\xi^*.
\]
Let $V_\Delta \in {\rm Hom}(id,\gamma)$, $W_\Delta \in {\rm Hom}(\gamma, 
\gamma^2)$ be isometries defined by 
\begin{eqnarray*}
 &&V_\Delta=V_{id_M}, \\
 &&W_\Delta=\sum_{\xi,\eta,\zeta \in \Delta_0} 
\sqrt{\frac{d(\xi)d(\eta)}{\dim \Delta d(\zeta)}} V_\xi 
\hat{\xi}(V_\eta) T_{\xi,\eta}^\zeta V_\zeta^*,
\end{eqnarray*}
where $T_{\xi,\eta}^\zeta=\sum_{i=1}^{N_{\xi,\eta}^\zeta} T(^\zeta_{\xi,\eta})_i
\otimes j(T(^\zeta_{\xi,\eta})_i)$.
Then, one can construct a subfactor $B_\Delta$ of $A$ such that 
$\gamma_\Delta: A \longrightarrow B_\Delta$ is the canonical 
endomorphism of the inclusion $A \supset B_\Delta$ \cite{LRe}. We call the 
inclusion $A \supset B_\Delta$ the {\it Longo-Rehren inclusion.} 

In a similar manner, we can construct the Longo-Rehren inclusion $A 
\supset B_{\hat{\Delta}}$. By their constructions, we have the inclusions 
$A \supset B_{\Delta} \supset B_{\hat{\Delta}}$. 

We define $D(\Delta)$ to be the set of endomorphisms $\rho \in 
{\rm End}(B_\Delta)_0$ such that $[\iota_\Delta][\rho]$ is a finite direct sum 
of sectors in the decompositions of $\{[\xi \otimes 
id^{op}][\iota_\Delta]\}_{\xi \in \Delta_0}$, where $\iota_\Delta$ is the 
inclusion map $\iota_\Delta:B_\Delta \hookrightarrow A$. We call 
$D(\Delta)$ the {\it quantum double} of $\Delta$. (For a categorical 
interpretation of the quantum double, see \cite{Muger1}.) In Corollary 7.2 
\cite{Izumi1}, it is proved that $D(\hat{\Delta})$ is equivalent to 
$\hat{\Delta} \otimes \hat{\Delta}^{op}$ as modular categories.

\begin{prop}\label{crossedproduct}
We assume that $\Delta^d \cong U(G)$, where $G$ is an 
abelian group. Then, there exists an outer action $\alpha$ of $G$ 
on $B_{\hat{\Delta}}$ and the subfactor $B_\Delta \supset B_{\hat{\Delta}}$ is 
isomorphic to $B_{\hat{\Delta}} \rtimes_\alpha G \supset  B_{\hat{\Delta}}$.
\end{prop}

\noindent
{\it Proof.} 
Let $\iota_1:B_{\hat{\Delta}} \hookrightarrow B_{\Delta}$ be the inclusion map. Then, by 
Theorem 7.4 in \cite{Izumi1}, we have $[\bar{\iota}_1 \iota_1]=\oplus_{\xi 
\in \Delta'_0} [\widehat{\rho_{\xi, \bar{\xi}}}^{+-}]$. By the 
minimality of the non-degenerate extension, this is
$[\bar{\iota}_1 \iota_1]=\oplus_{\xi \in \Delta^d_0} 
[\widehat{\rho_{\xi, \bar{\xi}}}^{+-}]$. Since $G \cong \hat{G}$ as groups and 
$d(\widehat{\rho_{\xi, \bar{\xi}}}^{+-})=d(\rho_\xi)=1$ for each $\xi \in 
\hat{\Delta}_0$, $\widehat{\rho_{\xi, \bar{\xi}}}^{+-}$ is an 
automorphism labeled by $G$. Then, by Theorem 4.1 in \cite{Iz}, there exists an 
outer action 
$\alpha$ of $G$ on $B_{\hat{\Delta}}$ and the dual inclusion of $B_\Delta \supset 
B_{\hat{\Delta}}$ is $B_\Delta \supset {B_{\Delta}}^G$. Hence, $B_\Delta \supset 
B_{\hat{\Delta}}$ is isomorphic to $B_{\hat{\Delta}} \rtimes_\alpha G \supset  
B_{\hat{\Delta}}$. 
\hfill $\Box$

\begin{thm}\label{main1}
 Let $D(\Delta)$ be the quantum double of $\Delta$. Then, under the assumptions 
in Proposition \ref{crossedproduct}, 
 $D(\Delta)=(\hat{\Delta} \otimes \hat{\Delta}^{op} \cap {\Delta^d}') 
 \rtimes \Delta^d$, where the embedding $\iota_{\Delta^d}:\Delta^d \hookrightarrow  \hat{\Delta} \otimes 
\hat{\Delta}^{op}$ is given by $\iota_{\Delta^d} (\sigma)=(\sigma, 
\sigma^{op})$. 
\end{thm}

\noindent
{\it Proof.} 
First, we may assume that $D(\hat{\Delta})=\hat{\Delta} \otimes 
\hat{\Delta}^{op}$ by thanks to Corollary 7.2 in \cite{Izumi1}. By its 
construction, $M \rtimes_\alpha G$ can be viewed as $M \rtimes \hat{G}$. Then, 
we may apply Proposition \ref{crossedproduct} to $\hat{\Delta} \otimes 
\hat{\Delta}^{op} \cap {\Delta^d}'$ to get the crossed product $(\hat{\Delta} \otimes 
\hat{\Delta}^{op} \cap {\Delta^d}')  \rtimes \Delta^d$ in 
${\rm End}(B_{\Delta})_0$. 

By Lemma 7.6 in \cite{Izumi1}, the image of the $\alpha$-induction in 
Proposition \ref{muger} is in $D(\Delta)$. Thus, $(\hat{\Delta} \otimes 
\hat{\Delta}^{op} \cap {\Delta^d}')  \rtimes \Delta^d$ is a full 
subcategory of $D(\Delta)$. 

We compute the dimension of $(\hat{\Delta} \otimes \hat{\Delta}^{op} \cap 
{\Delta^d}')  \rtimes \Delta^d$.
\begin{eqnarray*}
 \dim (\hat{\Delta} \otimes \hat{\Delta}^{op} \cap {\Delta^d}')  \rtimes \Delta^d
&=& \frac{\dim \hat{\Delta} \otimes \hat{\Delta}^{op} \cap {\Delta^d}'}{\dim \Delta^d}
= \frac{\dim \hat{\Delta} \otimes \hat{\Delta}^{op}}{(\dim \Delta^d)^2} \\
&=& \left(\frac{\dim \hat{\Delta}}{\dim \Delta^d}\right)^2 
= (\dim \Delta)^2 \\
&=& (\dim D(\Delta))^2,
\end{eqnarray*}
where we used the minimality of the extension $\hat{\Delta} \supset 
\Delta$ for the fourth equality.

Thus, $\dim D(\Delta)=\dim (\hat{\Delta} \otimes \hat{\Delta}^{op} \cap {\Delta^d}')  
\rtimes \Delta^d$, and this implies $D(\Delta)=(\hat{\Delta} \otimes \hat{\Delta}^{op} \cap {\Delta^d}') 
 \rtimes \Delta^d$. \hfill $\Box$
\section{Application to the Reshetikhin-Turaev invariants for 3-manifolds}
We apply Theorem \ref{main1} to the Reshetikhin-Turaev invariant of 
3-manifolds constructed from the quantum double $D(\Delta)$ to get a simpler 
description of it in this case. Before we state Theorem, we collect some general 
results on a premodular category.

\begin{lem}\label{lem1}
Let ${\cal M}$ be a premodular category, $\cal P$ the non-degenerate 
extension of $\cal M$ and $\cal D$ be degenerates of $\cal M$, i.e., 
${\cal D}={\cal M} \cap {\cal M}'$. Then, we have 
\begin{equation}
 \sum_{\omega \in {\cal M}_0} N_{\eta \bar{\zeta}}^\omega d(\omega) 
= d(\eta \bar{\zeta}) \chi_{{\cal  M}} (\eta \bar{\zeta}),
\end{equation}
where $\chi_{\cal M}(\xi)=1\ \text{if}\ \xi \in {\cal M}, 0\ \text{otherwise}$.
\end{lem}

\noindent
{\it Proof.}
We compute $\sum_{\xi \in {\cal D}_0} S'(\xi, \eta) S'(\xi, \bar{\zeta})$ in  
different ways. 

On one hand, we have 
\begin{eqnarray*}
 \sum_{\xi \in {\cal D}_0} S'(\xi, \eta) S'(\xi,\bar{\zeta}) 
&=& \sum_{\xi \in {\cal D}_0} \sum_{\omega \in {\cal P}_0} d(\xi) 
N_{\eta \bar{\zeta}}^\omega 
S'(\xi, \omega) \\
&=& \sum_{\omega \in {\cal P}_0} N_{\eta \bar{\zeta}}^\omega 
\sum_{\xi \in {\cal D}_0} d(\xi) 
S'(\xi, \omega)  \\
&=& \sum_{\omega \in ({{\cal P} \cap {\cal D}'})_0} N_{\eta \bar{\zeta}}^\omega 
d(\omega) \dim {\cal D} 
\end{eqnarray*}
where we used $\sum_{\xi \in {\cal D}_0} d(\xi) S'(\xi, \omega)=d(\omega) 
\dim {\cal D} \chi_{{\cal P} \cap {\cal D}'} (\omega)$ 
in Lemma 2.13 in \cite{Muger3} for the third equality. 

On the other hand,
\begin{eqnarray*}
\sum_{\xi \in {\cal D}_0} S'(\xi, \eta) S'(\xi, \bar{\eta}) 
&=& \sum_{\xi \in {\cal D}_0} S'(\xi, \eta \bar{\zeta}) d(\xi) \\
&=& d(\eta \bar{\zeta}) \dim {\cal D} \chi_{{\cal P} \cap {\cal D}'}
(\eta \bar{\zeta}). 
\end{eqnarray*}

Thus, $\sum_{\omega \in ({\cal P} \cap {\cal D}')_0} N_{\eta \bar{\zeta}}^\omega 
d(\omega) =d(\eta \bar{\zeta})  \chi_{{\cal P} \cap {\cal D}'}(\eta \bar{\zeta})$ 
with ${\cal P} \cap {\cal D}' = {\cal M}$ implies the claim. \hfill $\Box$

Let $\cal C$ be a premodular category. 
Let $L$ be a framed link with $n$ components in the 3-sphere. 
We denote the invariant of 
the colored framed link by $F_{{\cal C}}(L,\lambda)$, where 
$\lambda=(\lambda_1, \cdots, \lambda_n) \in {\cal C}_0^n$. Set 
\[
 \{ L\}_{\cal C}=\sum_{\lambda \in {\cal C}_0^n} \prod_{i=1}^n 
 d(\lambda_i) F_{\cal C}(L; \lambda).
\]
We may assume that a closed 3-manifold $M$ is obtained from surgery 
along the framed link $L$ in the 3-sphere $S^3$. Namely, $M=\partial 
W_L$, where $W_L$ is the 4-manifold obtained by gluing $n$ 2-handles 
to the 4-ball $B^4$ along $L \subset S^3=\partial B^4$. We denote the 
signature of $W_L$ by $\sigma(L)$.

Let $\cal C$ be a modular category and we set 
$\Delta_{\cal C}=\sum_{\xi \in {\cal C}_0} t_\xi^{-1} d(\xi)^2$ and 
$D_{\cal C}=(\dim {\cal C})^{1/2}$. 
The Reshetikhin-Turaev invariant $\tau_{\cal C}$ is defined by 
\[
 \tau_{\cal C}(M) =(\Delta_{\cal C})^{\sigma(L)} 
D_{\cal  C}^{-\sigma(L)-n-1} \{ L \}_{\cal C}.
\] 
See \cite{Turaev} for the details of the definition.
\begin{lem}\label{lem3}
 Let ${\cal C}$ be a premodular category with ${\cal C} \cap {\cal C}'
={\cal D}$ and $L$ be a framed link with $n$ components. Then, we have
\[
 \{ L \}_{\cal C} =(\dim {\cal D})^n \{ L \}_{{\cal C} \rtimes {\cal D}}.
\] 
\end{lem}

\noindent
{\it Proof.} 
This is immediate from Remarques 2.1 1) and Proposition 3.7 1) in 
\cite{Br}. \hfill $\Box$

We now go back in the case of braided $C^*$-tensor categories $\hat{\Delta}$ 
and $\Delta$ associated with subfactors. Recall that we have assumed the 
minimality of the non-degenerate extension $\hat{\Delta}\supset \Delta$.
For $\lambda, \mu \in \hat{\Delta}$, we put 
\[
 [\lambda, \mu]_\Delta=\frac{1}{\dim \hat{\Delta}} \sum_{\nu \in \Delta_0} 
N_{\lambda \bar{\mu}}^\nu d(\nu).
\]

\begin{thm}\label{main2} 
Let $M$  be a closed 3-manifold obtained from surgery along the framed 
link $L$ with $n$ components. 
Then, the Reshetikhin-Turaev invariant for $D(\Delta)$ is given by 
\[
 \tau_{D(\Delta)}(M)=\frac{1}{\dim \Delta} 
\sum_{\lambda, \mu \in  \hat{\Delta}_0^n} \prod_{i=1}^n 
[\lambda_i, \mu_i]_{\Delta} F_{\hat{\Delta}}(L;\lambda) \overline{F_{\hat{\Delta}}(L;\mu)}.
\]
\end{thm}

\noindent
{\it Proof.} 
Since $\Delta_{D(\Delta)}=D_{D(\Delta)}$, we have
\[
\tau_{D(\Delta)}(M)= \frac{1}{(\dim \Delta)^{n+1}} \sum_{\tilde{\xi} \in D(\Delta)_0^n} 
\prod_{i=1}^n d(\tilde{\xi_i}) F_{D(\Delta)}(L;\tilde{\xi}).
\]
Then, by Theorem \ref{main1} and Lemma \ref{lem3}, 
\begin{equation}\label{eqn1}
 \tau_{D(\Delta)}(M)= \frac{1}{(\dim \Delta)^{n+1} (\dim \Delta^d)^n} 
\sum_{\tilde{\zeta} \in (\hat{\Delta} \otimes \hat{\Delta}^{op} \cap {\Delta^d}')_0^n} 
\prod_{i=1}^n d(\tilde{\zeta_i}) F_{D(\Delta)}(L;\tilde{\zeta}).
\end{equation}
We note that  for $\tilde{\zeta} \in (\Delta \otimes \Delta^{op} \cap 
{\Delta^d}')_0$ there exist $\lambda, \mu \in \hat{\Delta}_0$ such that 
$\tilde{\zeta}=\lambda \otimes \mu^{op}$. With this and Lemma \ref{lem1}, 
we have 
$
 d(\tilde{\zeta})\chi_{\hat{\Delta} \cap {\Delta^d}'} (\lambda \bar{\mu})
=d(\lambda) d(\bar{\mu}) \chi_{\hat{\Delta} \cap {\Delta^d}'} (\lambda \bar{\mu})
=\dim \hat{\Delta}\; [\lambda, \mu]_\Delta 
$ 

Hence, the right hand side of (\ref{eqn1}) is  
\[
\frac{(\dim \hat{\Delta})^n}{(\dim \Delta)^{n+1} (\dim \Delta^d)^n} 
\sum_{\lambda, \mu \in {\hat{\Delta}}_0^n} 
\prod_{i=1}^n [\lambda_i, \mu_i]_{\Delta} 
F_{\hat{\Delta} \otimes \hat{\Delta}^{op}} (L; \lambda_1 \otimes \mu_1^{op}, 
\lambda_2 \otimes \mu_2^{op}, \cdots, \lambda_n \otimes \mu_n^{op}).
\] 

Since we have the equality $F_{\hat{\Delta} \otimes \hat{\Delta}^{op}} 
(L; \lambda_1 \otimes \mu_1^{op}, \lambda_2 \otimes \mu_2^{op}, \cdots, 
\lambda_n \otimes \mu_n^{op})=F_{\hat{\Delta}} (L; \lambda_1, \cdots, 
\lambda_n) \overline{F_{\hat{\Delta}} (L; \mu_1, \cdots, \mu_n)}$ and 
the minimality of the non-degenerate extension $\hat{\Delta} \supset 
\Delta$, we have 
\[
 \tau_{D(\Delta)}(M)=\frac{1}{\dim \Delta} 
\sum_{\lambda, \mu \in  \hat{\Delta}_0^n} \prod_{i=1}^n 
[\lambda_i, \mu_i]_{\Delta} F_{\hat{\Delta}}(L;\lambda) \overline{F_{\hat{\Delta}}(L;\mu)}.
\]
\hfill $\Box$

\begin{rem}
 With Theorem 5.2 in \cite{KSW}, which claims that 
the Turaev-Viro-Ocneanu invariant for $\Delta$ is equal to the Reshetikhin-Turaev 
invariant for $D(\Delta)$, Theorem \ref{main2} proves a slightly different statement 
of Theorem 3.2 in \cite{Ocneanu} in the special case, although Ocneanu claims that 
it also holds true for $G$, a non-abelian group. 
\end{rem}

\end{document}